\documentclass[10pt]{article}
\usepackage{dcolumn}
\usepackage{bm}
\usepackage{verbatim}       

\usepackage[dvips]{graphicx}

\usepackage{amsthm}
\usepackage{amssymb}

\usepackage{bm}
\usepackage{amstext}
\usepackage{amsmath}
\usepackage{amsfonts}
\usepackage{units}
\usepackage{mathrsfs}
\newcommand{\bd}{\begin{definition}}                
\newcommand{\ed}{\end{definition}}                  
\newcommand{\bc}{\begin{corollary}}                 
\newcommand{\ec}{\end{corollary}}                   
\newcommand{\bl}{\begin{lemma}}                     
\newcommand{\el}{\end{lemma}}                       
\newcommand{\bp}{\begin{proposition}}            
\newcommand{\ep}{\end{proposition}}                
\newcommand{\bere}{\begin{remark}}                  
\newcommand{\ere}{\end{remark}}                     

\newcommand{\bt}{\begin{theorem}}
\newcommand{\et}{\end{theorem}}

\newcommand{\be}{\begin{equation}}
\newcommand{\ee}{\end{equation}}

\newcommand{\bit}{\begin{itemize}}
\newcommand{\eit}{\end{itemize}}
\newtheorem{theorem}{Theorem}[section]
\newtheorem{corollary}[theorem]{Corollary}
\newtheorem{lemma}[theorem]{Lemma}
\newtheorem{proposition}[theorem]{Proposition}
\theoremstyle{definition}
\newtheorem{definition}[theorem]{Definition}
\theoremstyle{remark}
\newtheorem{remark}[theorem]{Remark}
\newtheorem{example}[theorem]{Example}

\begin{document}

\title{Normally preordered spaces and utilities\thanks{This  version differs from that published in Order
as it contains two proofs of theorem \ref{bhp}.}}


\author{E. Minguzzi\thanks{
Dipartimento di Matematica Applicata ``G. Sansone'', Universit\`a
degli Studi di Firenze, Via S. Marta 3,  I-50139 Firenze, Italy.
E-mail: ettore.minguzzi@unifi.it} }

\date{}

\maketitle

\begin{abstract}
\noindent In applications it is useful to know whether a topological
preordered space is normally preordered. It is proved that every
$k_\omega$-space equipped with a closed preorder is a normally
preordered space. Furthermore, it is proved that second countable
regularly preordered spaces are perfectly normally preordered and
admit a countable utility representation.
\end{abstract}

\section{Introduction}

In applications such as  dynamical systems \cite{akin93}, general
relativity \cite{minguzzi09c} or microeconomics \cite{bridges95} it
is useful to know if a topological preordered space, usually a
topological manifold, is normally preordered.\footnote{Domain theory
\cite{gierz03} has applications to computer science and is related
in a natural way to ordered topological spaces. In this field a
topological space equipped with a closed order is called {\em
pospace} and a normally ordered space is called {\em monotone normal
pospace}.} The preorder arises from the orbit dynamics of the
dynamical system; from the causal preorder of the spacetime
manifold; or from the preferences of the agent in microeconomics.
The condition of preorder normality can be regarded as just one
first step in order to prove that the space is quasi-uniformizable
or even quasi-pseudometrizable in such a way that it admits order
completions and order compactifications.

The case of a preorder is often as important as the case of an
order. Indeed, dynamical systems are especially interesting in the
presence of dynamical cycles, and, analogously, spacetimes are
particularly interesting in presence of causality violations. Also
the case of a preorder is the usual one considered in microeconomics
as an agent may be indifferent with respect to  two possibilities in
the space of alternatives (prospect space) (actually the agent can
even be unable to compare them, this possibility is called
indecisiveness or incomparability \cite{aumann62}).


It is well known that a topological space equipped with a closed
order is Hausdorff. The removal of the antisymmetry condition for
the order suggests to remove the Hausdorff condition for the
topology. Indeed, quite often in applications, the preorder is so
tightly linked with the topology that one has that two points which
are indistinguishable according to the preorder (i.e. $x\le y$ and
$y \le x$) are also indistinguishable according to the topology, so
that even imposing the $T_0$ property could be too strong. This fact
is not completely appreciated in the literature on topological
preordered spaces. Nachbin, in his foundational book
\cite{nachbin65}, uses at some crucial step the Hausdorff condition
implied by the closed order assumption \cite[Theor. 4, Chap.
I]{nachbin65}. In this work we shall remove altogether the Hausdorff
assumption on the topology and in fact even the $T_0$ assumption.
The idea is that the separability conditions for the topology should
preferably come from their preorder versions and  should not be
added to the assumptions.

So far the only result which allows us to infer that a topological
preordered space is normally preordered is Nachbin's theorem
\cite[Theor. 4, Chap. I]{nachbin65}, which states that a compact
 space equipped with a closed order is normally ordered. There are other results of
this type \cite[Theor. 4.9]{fletcher82} but they assume the totality
of the order. Our main objective is to prove a result that holds at
least for topological manifolds and in the preordered case so as to
be used in the mentioned applications. Indeed, we shall prove that
the $k_\omega$-spaces equipped with a closed preorder are normally
preordered. Since topological manifolds are second countable and
locally compact and these properties imply the $k_\omega$-space
property, the theorem will achieve  our goal.

\subsection{Topological preliminaries}


Since in this work we do not assume Hausdorffness,  it is necessary
to clarify that in our terminology a topological space is {\em
locally compact} if every point admits a compact neighborhood.





A topological space $E$ is a {\em k-space}  if $O\subset E$ is open
if and only if, for every compact set $K\subset E$, $O\cap K$ is
open in $K$. We remark that  we are using here the definition given
in \cite{willard70}, thus we do not include Hausdorffness in the
definition as done in \cite[Cor. 3.3.19]{engelking89}. Using our
definition of local compactness it is not difficult to prove that
every first countable or locally compact space is a $k$-space
(modifying slightly the proof in \cite[Theor. 43.9]{willard70}).




A related notion is that of {\em $k_{\omega}$-space} which can be
characterized through the following property \cite{franklin77}:
there is a countable sequence $K_i$  of compact sets such that
$\bigcup_{i=1}^{\infty}K_i= E$ and for every $O\subset E$, $O$ is
open if and only if $O\cap K_i$ is open in $K_i$ with the induced
topology (again, here $E$ is not required to be Hausdorff). The
sequence $K_i$ is called {\em admissible}.  By replacing $K_i \to
\bigcup_{j=1}^i K_j$ one checks that it is possible to assume $K_i
\subset K_{i+1}$; also the replacement $K_i \to K\cup K_i$ shows
that  in the admissible sequence $K_1$ can be chosen to be any
compact set. We have the chain of implications: compact
$\Rightarrow$ $k_\omega$-space $\Rightarrow$ $\sigma$-compact
$\Rightarrow$ Lindel\"of,
and the fact that local compactness makes the last three properties
coincide.

We shall be interested on the behavior of the $k_\omega$-space
condition under quotient maps. Remarkably, the next proof shows that
the Hausdorff condition in Morita's theorem \cite[Lemmas
1-4]{morita56} can be dropped.


\begin{theorem} \label{khg}
Every $\sigma$-compact locally compact (locally compact) space is a
$k_\omega$-space (resp. $k$-space). The quotient of a
$k_\omega$-space ($k$-space) is a $k_\omega$-space (resp.
$k$-space). Every $k_\omega$-space ($k$-space) is the quotient of a
$\sigma$-compact (resp. paracompact) locally compact space.
\end{theorem}

\begin{proof}
The first statement has been already mentioned.

Let ${K}_\alpha, \alpha \in \Omega$, be an admissible sequence
(resp. the family of all the compact sets) in $E$. Let $\pi: E\to
\tilde{E}$ be a quotient map and let
$\tilde{{K}}_\alpha=\pi({K}_\alpha)$; then the sets
$\tilde{{K}}_\alpha, \alpha \in \Omega$, form a countable family of
compact sets (resp. a subfamily of the family of all the compact
sets). Suppose that $C\subset \tilde{E}$ is such that $C\cap
\tilde{K}_\alpha$ is closed in $\tilde{K}_\alpha$. We have
$\pi^{-1}(C\cap \tilde{K}_\alpha)=\pi^{-1}(C)\cap
\pi^{-1}(\tilde{K}_\alpha)$; thus $\pi^{-1}(C)\cap
K_\alpha=\pi^{-1}(C)\cap \pi^{-1}(\tilde{K}_\alpha) \cap K_\alpha=
\pi^{-1}(C\cap \tilde{K}_\alpha)\cap K_\alpha$. Let $C_\alpha$ be a
closed set on $\tilde{E}$ such that $C_\alpha\cap \tilde{K}_\alpha=
C\cap \tilde{K}_\alpha$; then $\pi^{-1}(C)\cap K_\alpha=
\pi^{-1}(C_\alpha\cap \tilde{K}_\alpha)\cap
K_\alpha=\pi^{-1}(C_\alpha)\cap \pi^{-1}(\tilde{K}_\alpha)\cap
K_\alpha=\pi^{-1}(C_\alpha)\cap K_\alpha$. As the set on the
right-hand side is closed in $K_\alpha$ for every $\alpha$ we get
that $\pi^{-1}(C)$ is closed and hence $C$ is closed by the
definition of quotient topology. (Note that $\tilde{K}_\alpha$ is an
admissible sequence in the $k_\omega$-space case.)

 Let $K_\alpha$, $\alpha \in \Omega$, be an
admissible sequence (resp. the family of all the compact sets) in
$E$. Let $\tilde{K}_{\alpha}=\{(x,\alpha), x\in K_\alpha\}$,
$E'=\cup_{\alpha} \tilde{K}_\alpha$ be the disjoint union
\cite{willard70} and let $g:E' \to E$ be the map given by
$g((x,\alpha))=x$ so that $\varphi_{\alpha}:=g
\vert_{\tilde{K}_\alpha}: \tilde{K}_\alpha \to K_\alpha$ is a
homeomorphism. Let $C\subset E$ be such that $g^{-1}(C)$ is closed
then for each $\alpha$, $g^{-1}(C)\cap \tilde{K}_{\alpha}$ is closed
in $\tilde{K}_{\alpha}$ thus $g(g^{-1}(C)\cap
\tilde{K}_\alpha)=\varphi_\alpha(g^{-1}(C)\cap
\tilde{K}_\alpha)=C\cap K_\alpha$ is closed in $K_\alpha$ because
$\varphi_\alpha$ is a homeomorphism. Thus $C$ is closed and hence
$g$ is a quotient map. It is trivial to check that $E'$ is
$\sigma$-compact (resp. paracompact) and locally compact.
\end{proof}

\subsection{Order theoretical preliminaries}

For a topological preordered space $(E,\mathscr{T},\le)$ our
terminology and notation follow Nachbin \cite{nachbin65}. Thus by
{\em preorder} we mean a reflexive and transitive relation. A
preorder is an {\em order} if it is antisymmetric. With $i(x)=\{y:
x\le y\}$ and $d(x)=\{y: y\le x\}$ we denote the increasing and
decreasing hulls. A topological preordered space is {\em semiclosed
 preordered} if $i(x)$ and $d(x)$ are closed for every $x \in E$, and
it is {\em closed preordered} if the graph of the preorder
$G(\le)=\{(x,y): x\le y\}$ is closed. A subset $S\subset E$, is
called {\em increasing} if $i(S)=S$ and {\em decreasing} if
$d(S)=S$. $I(S)$ denotes the smallest closed increasing set
containing $S$, and $D(S)$ denotes the smallest closed decreasing
set containing $S$. It is understood that the set inclusion is
reflexive,  $X\subset X$.

A topological preordered space  is a {\em normally preordered space}
if it is semiclosed preordered and for every closed decreasing set
$A$ and closed increasing set $B$ which are disjoint, $A\cap
B=\emptyset$, it is possible to find an open decreasing set $U$ and
an open increasing set $V$ which separate them, namely $A\subset U$,
$B\subset V$, and $U\cap V=\emptyset$.

A {\em regularly preordered space} is a semiclosed preordered space
such that, if $x \notin B$ with $B$ a closed increasing set, there
is an open decreasing set $U\ni x$ and an open increasing set
$V\supset B$, such that $U\cap V=\emptyset$, and analogously the
dual property must hold for $y \notin A$ with $A$ a closed
decreasing set.

We have the implications: normally preordered space $\Rightarrow$
regularly preordered space $\Rightarrow$  closed preordered space
$\Rightarrow$ semiclosed preordered space.  A topological preordered
space is {\em convex} \cite{nachbin65} if for every $x\in E$, and
open set $O\ni x$, there are an open decreasing set $U$ and an open
increasing  set $V$ such that $x\in U\cap V\subset O$.


\section{Preorders on compact spaces and $k_\omega$-spaces}

We are interested in establishing in which way compactness and
countability assumptions improve the preorder separability
properties  of a topological preordered space. The example below
shows that these conditions do not promote semiclosed preordered
spaces to closed preordered spaces, not even under convexity, and
thus that the closed preorder property is in fact much more
interesting since, as we shall see, it allows us  to reach better
separability properties.


\begin{example}
Let $E=[0,1]^2$ with the usual product topology $\mathscr{T}$.
Evidently $(E,\mathscr{T})$ has very good topological properties: it
is second countable, Hausdorff, compact and even complete with
respect to the Euclidean metric. Let $(x,y)$ be coordinates on $E$
and let $\le$ be the order defined as follows
\begin{align*}
i((x,y))=\{(x',y'): \ x'=x\textrm{ and},& \textrm{ if } x>0, \ y\le
y';
\\&\textrm{ if } x=0 \textrm{ and } y\le 1/2, \ y\le y'\le 1/2;\\&
\textrm{ if } x=0 \textrm { and } y>1/2,\  y=y'\}.
\end{align*}
With this choice $\le$ is completely determined, and
$(E,\mathscr{T},\le)$ can be checked to be a convex semiclosed
ordered space which is not a closed ordered space.
\end{example}

We need to state the next two propositions that generalize to
preorders two corresponding propositions due to Nachbin \cite[Prop.
4,5, Chap. I]{nachbin65}. Actually the proofs given by Nachbin for
the case of an order work also in this case without any
modification. For this reason they are omitted.

\begin{proposition} \label{pr1}
Let $E$ be a  closed preordered space. For every compact $K\subset
E$, we have $d(K)=D(K)$ and $i(K)=I(K)$, that is, the decreasing and
increasing hulls are closed.
\end{proposition}

\begin{proposition} \label{pr2}
Let $E$ be a compact closed preordered  space. Let $F\subset V$
where $F$ is increasing and $V$ is open, then there is an open
increasing set $W$ such that $F\subset W\subset V$. An analogous
statement holds in the decreasing case.
\end{proposition}

Every compact space equipped with a closed order is normally ordered
\cite[Theor. 4, Chap. I]{nachbin65}. We shall need a slightly
stronger statement.

\begin{theorem} \label{nqs2}
Every compact space $E$ equipped with a closed preorder  is a
normally preordered space.
\end{theorem}

\begin{proof}
If $A\cap B=\emptyset$ with $A$ closed decreasing and $B$ closed
increasing, for every $x\in A$ and $y\in B$ we have $d(x)\cap
i(y)=\emptyset$, thus there is (see \cite[Prop. 1, Chap.
I]{nachbin65}) a decreasing neighborhood $U(x,y)$ of $x$ and an
increasing neighborhood $V(y,x)$ of $y$ (they are not necessarily
open) such that $U(x,y)\cap V(y,x)=\emptyset$. Since $A$ and $B$ are
closed subsets of a compact set they are compact and $\cup_{y\in B}
V(y,x)\supset B$ thus there are points $y_i\in B$, $i=1,\cdots, k$,
such that, defined the increasing neighborhood of $B$, $V(x):=\cup_i
V(y_i,x)$, and the decreasing neighborhood of $x$, $U(x):=\cap_i
U(x,y_i)$ we have $U(x)\cap V(x)=\emptyset$. The neighborhoods
$U(x)$, $x\in A$ are such that $\cup_{x\in A} U(x)\supset A$ thus we
can find $x_j$, $j=1,\cdots, n$ such that defined the decreasing
neighborhood of $A$, $U':=\cup_j U(x_j)$, and the increasing
neighborhood of $B$, $V':=\cap_j V(x_j)$, we have $U'\cap
V'=\emptyset$.  Finally, by Prop. \ref{pr2} there are an open
decreasing set $U$ such that $A\subset U\subset U'$, and an open
increasing set $V$ such that $V'\supset V\supset B$, from which the
thesis i.e. $U\cap V=\emptyset$.
 \end{proof}
%
%

A subset $S\subset E$ with the induced topology $\mathscr{T}_S$ and
the induced preorder $\le_S$ is a topological preordered space hence
called {\em subspace}. In general it is not true that every open
increasing (decreasing) set on $S$ is the intersection of an  open
increasing (resp. decreasing) set on $E$ with $S$. If this is the
case $S$ is called a {\em preordered subspace}
\cite{priestley72,mccartan68,kunzi90}.

\begin{proposition} \label{bog}
Every subspace of a (semi)closed preordered space is a
(semi)\-closed preordered space.
\end{proposition}

\begin{proof}

Let $(E,\mathscr{T},\le)$ be semiclosed preordered. Given $x\in S$,
let $i_S(x)=\{y\in S: x \le_S y\}=\{y\in S: x \le y\}$. From this
expression we have $i_S(x)=i(x)\cap S$, thus $i_S(x)$ is closed in
the induced topology. Analogously, $d_S(x)$ is closed in the induced
topology, that is $(S,\mathscr{T}_S,\le_S)$ is semiclosed
preordered.

Let $(E,\mathscr{T},\le)$ be closed preordered, thus $G(\le)$ is
closed, and hence the graph of $\le_S$, $G\cap (S\times S)$ is
closed in the induced (product)  topology $(\mathscr{T}\times
\mathscr{T})_{S\times S}$. The equality $(\mathscr{T}\times
\mathscr{T})_{S\times S}=(\mathscr{T}_S\times \mathscr{T}_S)$ proves
that the graph of $\le_S$ is closed.
\end{proof}

\begin{proposition} \label{bax}
In a closed preordered space every compact subspace $S$ is a
preordered subspace.
\end{proposition}

\begin{proof}
If $A \subset S$ is closed decreasing then it is compact, thus
$d(A)$ is closed decreasing and such that $A=d(A)\cap S$. The proof
in the increasing case is analogous.
\end{proof}

Using the previous results it is possible to follow the strategy of
the proof   that a Hausdorff  $k_\omega$-space is a normal space
\cite{franklin77} in order to obtain the next theorem.

\begin{theorem} \label{bhp}
Let $(E,\mathscr{T},\le)$ be a $k_\omega$-space endowed with a
closed preorder,\footnote{We shall call these spaces: closed
preordered $k_\omega$-spaces. However, this terminology is different
from that used in \cite{kopperman05} where in their closed ordered
$k$-spaces the term ``$k$-space''  includes an additional condition
on the upper and lower topologies.} then $(E,\mathscr{T},\le)$ is a
normally preordered space.
\end{theorem}

We shall give another proof in the next section. It must be noted
that Hausdorff locally compact spaces need not be normal thus in
theorem \ref{bhp} the $k_\omega$-space condition cannot be weakened
to the $k$-space condition (consider the discrete order).

\begin{proof}$\empty$\!\!\!\footnote{This first proof does not appear in the version of the paper published in Order. }
Let $K_n$, $K_n\subset K_{n+1}$, $\bigcup_n K_n=E$, be an admissible
sequence of compact sets (i.e. that appearing in the definition of
$k_\omega$-space). A set $O$ is open if and only if $O\cap K_n$ is
open in $K_n$.

Let $A,B$ be respectively a closed decreasing and a closed
increasing set such that $A\cap B=\emptyset$ and denote $A_n=A\cap
K_n$, $B_n=B\cap K_n$. By proposition \ref{bog} and theorem
\ref{nqs2} every $K_n$ with the induced topology and preorder is a
normally preordered space. Let $\tilde{A}_1=A_1$, $\tilde{B}_1=B_1$
and let $U_1,V_1\subset K_1$ be respectively open decreasing and
open increasing sets in $K_1$ (with respect to the induced topology
and preorder) such that $D_1(U_1)\cap I_1(V_1)=\emptyset$,
$\tilde{A}_1\subset U_1$, $\tilde{B}_1\subset V_1$, where $D_1, I_1$
are the closure operators for the preordered space $K_1$. They
exist, it suffices to apply the preordered normality of the space
three times. The set $I_1(V_1)$ being a closed set in $K_1$ is,
regarded as a subset of $E$, the intersection between a closed set
and a compact set thus it is compact. The set $i(I_1(V_1))$ is
closed (Prop. \ref{pr1}) and analogously, $d(D_1(U_1))$ is closed.
Now, let us consider the closed increasing set on $K_2$ given by
$\tilde{B}_2=[i(I_1(V_1))\cap K_2]\cup B_2$ and the closed
decreasing set given by $\tilde{A}_2=[d(D_1(U_1))\cap K_2]\cup A_2$
(see figure \ref{bhp}).

\begin{figure}[ht]
\centering
\includegraphics[width=7cm]{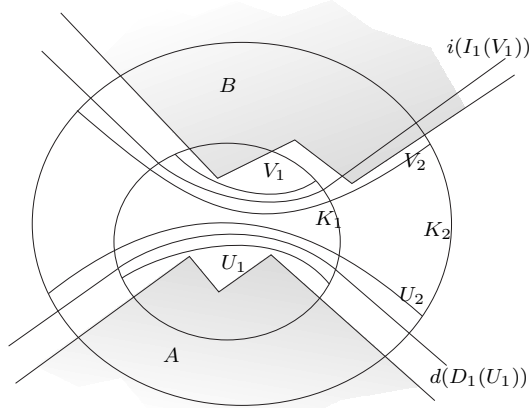}
\caption{The idea of the proof of theorem \ref{bhp}.} \label{pro}
\end{figure}

They are disjoint because $i(I_1(V_1))\cap d(D_1(U_1))=\emptyset$
(otherwise $I_1(V_1)\cap D_1(U_1)\ne\emptyset$, a contradiction),
$B_2\cap A_2=\emptyset$, and $B_2\cap d(D_1(U_1))=\emptyset$ as
$i(B_2)\cap D_1(U_1)\subset B_2\cap D_1(U_1)\subset V_1\cap D_1(U_1)
=\emptyset$. Analogously, $i(I_1(V_1))\cap A_2=\emptyset$. Thus,
arguing as before we can find $U_2,V_2\subset K_2$ respectively open
decreasing and open increasing sets in $K_2$ (with respect to the
induced topology and preorder) such that $D_2(U_2)\cap
I_2(V_2)=\emptyset$, $\tilde{A}_2\subset U_2$, $\tilde{B}_2\subset
V_2$. Continuing in this way we  define at each step
$\tilde{B}_{j+1}=[i(I_j(V_j))\cap K_{j+1}]\cup B_{j+1}$,
$\tilde{A}_{j+1}=[d(D_j(U_j))\cap K_{j+1}]\cup A_{j+1}$. Arguing as
before $\tilde{A}_{j+1}, \tilde{B}_{j+1}$ are disjoint closed
decreasing and closed increasing subsets of $K_{j+1}$ and since the
latter is a normally preordered space there are
$U_{j+1},V_{j+1}\subset K_{j+1}$ respectively open decreasing and
open increasing sets in $K_{j+1}$ such that $D_{j+1}(U_{j+1})\cap
I_{j+1}(V_{j+1})=\emptyset$, $\tilde{A}_{j+1}\subset U_{j+1}$,
$\tilde{B}_{j+1}\subset V_{j+1}$.

Note that $V_j \subset \tilde{B}_{j+1} \subset V_{j+1}$ and
analogously, $U_j\subset U_{j+1}$. Let $V=\bigcup_j V_j$ and
$U=\bigcup_j U_j$. The set $V$ contains $B$ because $B_j\subset
\tilde{B}_j\subset V_j$ thus $B=\bigcup_j B_j\subset V$.
Analogously, $U$ contains $A$.

The set $V$ is open because $V\cap K_s=\bigcup_{j\ge 1} (V_j\cap
K_s)=\bigcup_{j\ge s} (V_j\cap K_s)$, and the set $V_j\subset K_j$
is open in $K_j$ so that, since for $j\ge s$, $K_s \subset K_j$,
$V_j\cap K_s$ is open in $K_s$ and so is the union $V\cap K_s$. The
 $k_\omega$-space property implies that $V$ is open. Analogously,
$U$ is open.

Finally, let us prove that $V$ is increasing. Let $x\in V$ then
there is some $j\ge 1$ such that $x\in V_j \subset K_j$. Let $y\in
i(x)$, then we  can find some $r\ge j$ such that $y \in K_r$. Since
$V_j \subset V_r$, $x\in V_r$, and since $V_r$ is increasing on
$K_r$, $y \in V_r$ thus $y\in V$. Analogously, $U$ is decreasing
which completes the proof.

\end{proof}

\begin{corollary} \label{nai}
Every  locally compact $\sigma$-compact space equipped with a closed
preorder is a normally preordered space.
\end{corollary}

We can regard as a corollary  the known result \cite{franklin77},

\begin{corollary}
Every Hausdorff  $k_\omega$-space is a normal space.
\end{corollary}

\begin{proof}
Apply theorem \ref{bhp} to the discrete order and use the fact that
the Hausdorff condition is equivalent to the closure of the graph of
the discrete order, namely the diagonal $\Delta=\{(x,y): x=y\}$.
\end{proof}


By a result due to Milnor \cite[Lemma 2.1]{milnor56} Hausdorff
$k_\omega$-spaces are finitely productive. The locally compact
$\sigma$-compact spaces are finitely productive even if they are not
Hausdorff \cite{bourbaki66}, thus closed preordered locally compact
$\sigma$-compact  spaces are finitely productive. The closed
preordered locally compact $\sigma$-compact spaces provide our main
example because, by using the non-Hausdorff generalization of
Morita's theorem, we get  the following result.

\begin{corollary}
Every closed preordered $k_\omega$-space ($k$-space) $E$ is the
quotient of a closed preordered $\sigma$-compact (resp. paracompact)
locally compact space $E'$, $\pi: E'\to E$, in such a way that
$G(\le')=(\pi\times \pi)^{-1}(G(\le))$.
\end{corollary}

\begin{proof}
Obvious from theorem \ref{khg} because $\le'$ so defined is a closed
preorder.
\end{proof}

\section{A proof based on an extension theorem}
A  function $f: E \to \mathbb{R}$  is {\em isotone}, if $x\le y
\Rightarrow f(x)\le f(y)$. Nachbin proved that a continuous isotone
function $f: S \to [0,1]$ defined on a compact subset $S$ of a
normally ordered space $E$ can be extended to a function $F: E \to
[0,1]$ on the whole space preserving continuity and the isotone
property \cite[Theor. 6, Chap. I]{nachbin65}.\footnote{The fact that
the theorem holds with the functions $f,F,$ taking values in $[0,1]$
is evident from Nachbin's proof but is not stated in the original
theorem.} Unfortunately, he uses the order condition (and hence the
implied Hausdorff condition) and we need therefore to generalize the
theorem to the preordered case. Levin gives a similar result for
closed preorders on a compact set $E$ \cite[Lemma 2]{levin83},
\cite{levin83b} \cite[Theor. 6.1]{levin79}, in the context of the
mass transfer problem, nevertheless we prefer to give a proof closer
in spirit to Nachbin's topological approach.

\begin{theorem} \label{glg}
Let $E$ be a normally preordered space and let $S$ be a subspace.
Let $f:S\to [0,1]$ be continuous and isotone on $S$. In order that
the function $f$ be extendible to a continuous isotone function $F:
E\to [0,1]$, it is necessary and sufficient that
\begin{equation} \label{mjb}
 \xi, \xi' \in [0,1], \  \xi<\xi' \Rightarrow
D(f^{-1}([0,\xi]))\cap I(f^{-1}([\xi',1]))=\emptyset.
\end{equation}
\end{theorem}

The proof of this theorem is the same as the proof of \cite[Theor.
2]{nachbin65}. Indeed, in the body of that proof the image  of the
original and extended functions is in [0,1] and, rather
surprisingly, a close check of the proof (also of the omitted
continuity part \cite[Chap. 4, Lemma 3]{kelley55}) shows that it
does not depend on the closure condition on $S$ which is imposed in
the theorem statement.

\begin{remark}
Instead of rechecking  Nachbin's proof \cite[Theor. 2]{nachbin65}
one might just follow the argument below to show that $f$ in theorem
\ref{glg} can be extended to $\overline{S}$ preserving continuity,
the isotone property and Eq. (\ref{mjb}). Then one could apply
Nachbin's result \cite[Theor. 2]{nachbin65} and get theorem
\ref{glg}.

Under the hypotheses of Theorem \ref{glg}, let $x\in
\overline{S}\backslash S$. Let $(U_\alpha)_\alpha$ be a neighborhood
basis of $x$. Every $U_\alpha$ meets $S$. Thus $(U_\alpha \cap
S)_\alpha$ is a filter basis on $S$. Let $\eta=\liminf_\alpha
f(U_\alpha \cap S)$ and $\eta' = \limsup_\alpha f(U_\alpha \cap S)$.
If $\eta< \eta'$ we can find $\eta<\xi<\xi'<\eta'$. By hypothesis
$D(f^{-1}([0, \xi])) \cap I(f^{-1}([\xi', 1]))=\emptyset$. But $x$
belongs to this intersection as it belongs to $\overline{f^{-1}([0,
\xi])} \cap \overline{f^{-1}([\xi', 1])}$, which gives a
contradiction. Thus $\eta=\eta'$ and $f$ can be extended
continuously to $x$. The extended function $\tilde{f}:
\overline{S}\to [0,1]$ is isotone. Indeed, the inequality
$\tilde{f}(y)<\tilde{f}(x)$ implies that there are $\zeta<\zeta'$
such that $\tilde{f}(y)<\zeta<\zeta'<\tilde{f}(x)$. By continuity of
$\tilde{f}$, $x\in \overline{f^{-1}([\zeta', 1])}$ and $y\in
\overline{f^{-1}([0, \zeta])}$ and using $D(f^{-1}([0, \zeta])) \cap
I(f^{-1}([\zeta', 1]))=\emptyset$, we conclude that $x\nleq y$.
Finally, let $\xi, \xi' \in [0,1]$,  $\xi<\xi'$ and choose $\zeta,
\zeta' \in [0,1]$,  $\xi<\zeta<\zeta'<\xi'$ then
$D(f^{-1}([0,\zeta]))\cap I(f^{-1}([\zeta',1]))=\emptyset$ but $
D(\tilde{f}^{-1}([0,\xi])) \subset D(f^{-1}([0,\zeta]))$ and
$I(\tilde{f}^{-1}([\xi',1])) \subset I({f}^{-1}([\zeta',1]))$ from
which it follows $D(\tilde{f}^{-1}([0,\xi]))\cap
I(\tilde{f}^{-1}([\xi',1]))=\emptyset$.

\end{remark}

This generalization allows us to remove the closure condition in
\cite[Theor. 3]{nachbin65}, which therefore reads

\begin{theorem} \label{qin}
Let $E$ be a normally preordered space and let $S$ be a subspace
with the property that if $X,Y\subset S$ satisfy $D_S(X)\cap
I_S(Y)=\emptyset$ then $D(X)\cap I(Y)=\emptyset$. Then every
continuous isotone function $f:S\to [0,1]$ can be extended to a
continuous isotone function $F:E\to [0,1]$.
\end{theorem}

If we are in the discrete order case and $S$ is closed, the
assumption on theorem \ref{qin} is satisfied, thus one recovers
Tietze's extension theorem for bounded functions. We can now prove
(note that $S$ need not be closed; see also  remark \ref{mpa} for a
different proof).

\begin{theorem} \label{exl}
Let $E$ be a normally preordered  space and let $S$ be a compact
subspace, then any continuous isotone function  $f: S\to [0,1]$  can
be extended to a continuous isotone function $F: E\to [0,1]$.
\end{theorem}

\begin{proof}
Let $X,Y\subset S$ be such that $D_S(X)\cap I_S(Y)=\emptyset$. The
set $D_S(X)$ is a closed subset of the compact set $S$ thus it is
compact which implies that, by Prop. \ref{pr1}, $d(D_S(X))$ is
closed. Analogously, $i(I_S(Y))$ is closed. Moreover, $d(D_S(X))\cap
i(I_S(Y))=\emptyset$ for if $y\in d(D_S(X))\cap i(I_S(Y))$ there are
$z\in D_S(X)$ and $x\in I_S(Y)$, such that $x\le y\le z$ which
implies $x\le z$ or $\emptyset \ne i_S(I_S(Y))\cap D_S(X)=I_S(Y)\cap
D_S(X)$, a contradiction. As a consequence, $D(X)\cap I(Y)\subset
d(D_S(X))\cap i(I_S(Y))=\emptyset$. The desired conclusion follows
now from theorem \ref{qin}.
\end{proof}

%
%
%
%
%
%
%

\begin{lemma}\label{vad} (extension and separation lemma)
Let $(E,\mathscr{T},\le)$ be a   closed preordered compact space.
Let $K$ be a (possibly empty) compact subset of $E$, let $A\subset E
$ be a closed decreasing set and let $B\subset E$ be a closed
increasing set such that $A\cap B=\emptyset$. Let $f:K \to [0,1]$ be
a continuous isotone function on $K$ such that $A\cap K\subset
f^{-1}(0)$ and $B\cap K \subset f^{-1}(1)$. Then there is a
continuous isotone function $F:E \to [0,1]$ which extends $f$ such
that $A\subset F^{-1}(0)$ and $B \subset F^{-1}(1)$.
\end{lemma}

\begin{proof}

By theorem \ref{nqs2} the topological preordered space $E$ is a
normally preordered space.  Since $A$ and $B$ are closed and hence
compact subsets of $E$, the set $K'=A\cup K\cup B$ is a compact
subset of $E$.  The function $f':K' \to [0,1]$ defined by
$f'\vert_{A}=0$, $f'\vert_{K}=f$, $f'\vert_{B}=1$, is isotone.

Let us prove that $f'$ is continuous on $K'$ with the induced
topology. Clearly $f'^{-1}([0,1])$ is closed in $K'$ as it coincides
with $K'$. We need only to prove that $f'^{-1}([\alpha,1])$, $\alpha
>0$, is closed in $K'$, the proof for the case $f'^{-1}([0,\beta])$,
$\beta <1$, being analogous. The set $f^{-1}([\alpha,1])$ being a
closed subset of $K$ is a compact subset of $E$ thus
$I(f^{-1}([\alpha,1]))=i(f^{-1}([\alpha,1]))$. But $A\cap
f^{-1}([\alpha,1])=\emptyset$ and $A$ is decreasing thus $A\cap
I(f^{-1}([\alpha,1]))=\emptyset$. Since $f$ is continuous on $K$
there is some closed set $C$ in $E$ such that
$f^{-1}([\alpha,1])=C\cap K$. The closed set $C'=C\cap
I(f^{-1}([\alpha,1]))$ has again the property
$f^{-1}([\alpha,1])=C'\cap K$ and is disjoint from $A$. Now we can
write
\[
f'^{-1}([\alpha,1])=B\cup (C'\cap K)=B\cup (C'\cap K'),
\]
which proves that $f'^{-1}([\alpha,1])$ is the union of two closed
subsets of $K'$. We conclude that $f'$ is continuous on $K'$. By
theorem \ref{exl} $f'$ can be extended to a continuous isotone
function $F: E\to [0,1]$, which is the desired function.
\end{proof}

\begin{theorem} (improved extension and separation result) \label{jqx}
Let $(E,\mathscr{T},\le)$ be a $k_\omega$-space equipped with a
closed preorder. Let $K$ be a (possibly empty) compact subset, let
$D$ be closed decreasing and let $I$ be closed increasing, $D\cap
I=\emptyset$. Let $f:K \to [0,1]$ be a continuous isotone function
on $K$  such that $D\cap K\subset f^{-1}(0)$ and $I\cap K \subset
f^{-1}(1)$. Then there is a continuous isotone function $F:E \to
[0,1]$ which extends $f$ such that $D\subset F^{-1}(0)$ and $I
\subset F^{-1}(1)$.
\end{theorem}

\begin{proof}
%


Let $K_i$,  $E=\bigcup_i K_i$,  be an admissible sequence according
to the definition of $k_\omega$-space. Without loss of generality we
can assume $K_i \subset K_{i+1}$ and $K_1=K$.

By theorem \ref{nqs2} and proposition \ref{bog} each subset $K_i$
endowed with the induced topology and preorder is a normally
preordered space. Define $f_1: K_1 \to [0,1]$,  by $f_1=f$.
We make the inductive assumption that there is a continuous isotone
function $f_i: K_i \to [0,1]$ such that $D\cap K_i\subset
f_i^{-1}(0)$ and $I\cap K_i \subset f_i^{-1}(1)$. Applying lemma
\ref{vad} to $E=K_{i+1}$ with the induced order,  $A=D\cap K_{i+1}$,
$B=I\cap K_{i+1}$, and compact subspace $K_i$, we get that there is
a continuous isotone function $f_{i+1}: K_{i+1} \to [0,1]$ which
extends $f_i$ such that $D\cap K_{i+1}\subset f_{i+1}^{-1}(0)$ and
$I\cap K_{i+1} \subset f_{i+1}^{-1}(1)$. We conclude that there is
an isotone function $F: E\to [0,1]$ defined by $F\vert_{K_i}=f_i$
such that $D\subset F^{-1}(0)$ and $I \subset F^{-1}(1)$.

We recall that in a $k_\omega$-space with admissible sequence $K_i$,
$K_i \subset K_{i+1}$, a function $g:E \to \mathbb{R}$ is continuous
if and only if for every $i$, $g\vert_{K_i}$ is continuous in the
subspace $K_i$. By construction, the function $f_i$ is continuous in
$K_i$ thus  $F$ is continuous.
\end{proof}

%
%

We  can  give a second proof to theorem \ref{bhp}.


\begin{proof}[Second Proof]
Let $D\subset E$ be a closed decreasing subset and let $I\subset E$
be a closed increasing subset such that $D\cap I=\emptyset$. Let
$K=\emptyset$, then by theorem \ref{jqx} there is a continuous
isotone function $F:E \to [0,1]$  such that $D\subset F^{-1}(0)$ and
$I \subset F^{-1}(1)$. The open sets $\{x: F(x)>1/2\}$ and $\{x:
F(x)<1/2\}$, are respectively open increasing, open decreasing,
disjoint and containing respectively $I$ and $D$, thus
$(E,\mathscr{T},\le)$ is a normally preordered space.
\end{proof}

\begin{example}
We give an example of normally preordered space which admits a
non-closed subset $S$, which  satisfies the assumptions of theorem
\ref{qin}. Let $E=\mathbb{R}\times S^1$, $S^1=[0,2\pi)$, be equipped
with the product preorder $\le$,  where $\mathbb{R}$ is endowed with
the usual order $\preceq$, and $S^1$ is given the indiscrete
preorder. Let the topology $\mathscr{T}$ on $E$ be the coarsest
topology which makes the projection on the first factor, $\pi: E\to
\mathbb{R}$, continuous. By theorem \ref{bhp} $E$ is normally
preordered (or use the remark \ref{dks} below, and the fact that
$(\mathbb{R},\preceq)$ is normally ordered). The subset
$S=[0,\pi]^2$ is compact, thus it satisfies the assumptions of
theorem \ref{qin}, but non-closed, its closure being
$\overline{S}=[0,\pi]\times S^1$. The subset $S=(0,\pi)^2$ is
non-closed and non-compact but it still satisfies the assumptions of
theorem \ref{qin}.
\end{example}

\section{The ordered quotient space}

Let us introduce the equivalence relation $x\sim y$ on $E$, given by
``$x\le y$ and $y \le x$''. Let $E/\!\!\sim$ be the quotient space,
$\mathscr{T}/\!\!\sim$ the quotient topology, and let $\lesssim$ be
defined by, $[x]\lesssim [y]$ if $x\le y$ for some representatives
(with some abuse of notation we shall denote with $[x]$ both a
subset of $E$ and a point on $E/\!\!\sim$). The quotient preorder is
by construction an order. The triple
$(E/\!\!\sim,\mathscr{T}/\!\!\sim,\lesssim)$ is a topological
ordered space and $\pi: E\to E/\!\!\sim$ is the continuous quotient
projection.

\begin{remark} \label{dks}
Taking into account the definition of quotient topology we have that
every open (closed) increasing (decreasing) set on $E$ projects to
an open (resp. closed) increasing (resp. decreasing) set on
$E/\!\!\sim$ and all the latter sets can be regarded as such
projections. As a consequence, $(E,\mathscr{T},\le)$ is a normally
preordered space (semiclosed preordered space, regularly preordered
space) if and only if $(E/\!\!\sim,\mathscr{T}/\!\!\sim,\lesssim)$
is a normally ordered space (resp. semiclosed ordered space,
regularly ordered space).
\end{remark}

\begin{remark} \label{mpa}
An alternative proof of theorem \ref{exl} uses the fact that
$E/\!\!\sim$ is normally ordered and $\pi(S)$ is compact, so that
$f$ can be passed to the quotient, extended using \cite[Theor. 6,
Chap. I]{nachbin65} and then lifted to $E$.
\end{remark}

The closed preordered property does not pass smoothly to the
quotient, but in the compact case and in the $k_\omega$-space case,
using theorem \ref{bhp} (\ref{nqs2}) and theorem \ref{khg} we
obtain.

\begin{corollary} \label{mbt}
If $E$ is a  closed preordered compact space then $E/\!\!\sim$ is a
 closed ordered compact space. If $E$ is a closed preordered
$k_\omega$-space then $E/\!\!\sim$ is a closed ordered
$k_\omega$-space.
\end{corollary}

\begin{proof}
The first statement is a trivial consequence of theorem \ref{nqs2}.
As for the second statement, by Theor. \ref{khg}
$(E/\!\!\sim,\mathscr{T}/\!\!\sim)$ is a $k_\omega$-space. Since $E$
is a closed preordered $k_\omega$-space then it is normally
preordered from which it follows that $E/\!\!\sim$ is a normally
ordered space and hence a closed ordered space.
\end{proof}

\begin{remark}
The first statement in the previous result is contained in
\cite[Lemma 1]{candeal95} but the proof is incorrect again for the
tricky Hausdorff condition which they inadvertently use   in the
quotient. Indeed, they argument as follows: they take a net
$([a_\alpha],[b_\alpha])\in G(\lesssim)$ converging to $([a],[b])$
and prove that a subnet $(a_\beta,b_\beta)\in G(\le)$ converges to
some pair $(a',b')\in G(\le)$. Since $\pi$ is continuous
$([a_\beta],[b_\beta])$ converges to $([a'],[b']) \in G(\lesssim)$
(and also to $([a],[b])$) but this does not mean that
$([a'],[b'])=([a],[b])$ as the uniqueness of the limit requires the
Hausdorff condition and this is assured only {\em after} it is
proved that  $E/\!\!\sim$ is a closed ordered space.
\end{remark}

\begin{remark}
In the Hausdorff case, the second statement in corollary \ref{mbt}
can be proved using the strategy contained in \cite{lawson71}. If
$E$ is a Hausdorff $k_\omega$-space then $E/\!\!\sim$ is a Hausdorff
$k_\omega$-space \cite[Prop. 2.3b]{lawson71}, then $\pi \times \pi$
is a quotient map \cite[Prop. 2.3a, 2.2]{lawson71} which implies
since $G(\le)=(\pi\times \pi)^{-1}(G(\lesssim))$ that $G(\lesssim)$
is closed. One can then work out the proof of theorem \ref{bhp} in
the ordered framework of $E/\!\!\sim$ using remark \ref{dks}.
\end{remark}

\section{The existence of continuous utilities}

Let  us write $x<y$ if $x\le y$ and $ y \not\le x$. A {\em utility}
is a function $f: E\to \mathbb{R}$ such that ``$x\sim y \Rightarrow
f(x)=f(y)$ and $x< y \Rightarrow f(x)<f(y)$''. We say that the
preorder admits a {\em representation} by a family of functions
$\mathcal{F}$ if ``$x\le y$ iff $\forall f\in \mathcal{F}, f(x)\le
f(y)$''. It is easy to prove that the preorder of a normally
preordered space is represented by the family of continuous isotone
functions \cite[Theor. 1]{nachbin65}. It is interesting to
investigate under which conditions a continuous utility or, more
strongly, a representation through continuous utilities exists. This
problem has been thoroughly investigated, especially in the
economics literature. The simplest approach passes through the
assumption that the preordered space under consideration is normally
preordered \cite{mehta88}, although alternative strategies based on
weaker hypothesis have also been investigated \cite{herden89}. The
representation of relations through isotone and utility functions is
still an active field of research \cite{ok02,bosi06,bosi10}.

%
%

The reader must be warned that some authors call {\em utility} what
we call {\em isotone} function \cite{herden02,bosi06,evren09} and
use the word {\em representation} for the existence of just one
utility, although one utility does not allow us to recover the
preorder. This unfortunate circumstance comes from the fact that in
economics most terminology was introduced in connection with the
total preorder case, that is, before the importance of the general
preorder case was recognized.

Let us recall that a {\em perfectly normally preordered} space is a
semiclosed preordered space such that if $A$ is a closed decreasing
set and $B$ is a closed increasing set with $A\cap B=\emptyset$ then
there is a continuous isotone function $f: E \to [0,1]$ such that
$A=f^{-1}(0)$ and $B=f^{-1}(1)$.  Clearly,  perfectly normally
preordered spaces are normally preordered spaces.

A closed decreasing (increasing) set $S$ is {\em
functionally-preordered closed} if there is a continuous isotone
function $f:E \to [0,1]$ such that $S=f^{-1}(0)$ (resp.
$S=f^{-1}(1))$. They can also be called decreasing (increasing) {\em
zero sets} as done in \cite{mccallion72}. A pair $(A,B)$, $A\cap
B=\emptyset$, is {\em functionally-preordered closed} if there is a
continuous isotone function $f:E \to [0,1]$ such that $A=f^{-1}(0)$
and $B=f^{-1}(1)$. The next result is stated without proof in
\cite{mccallion72}.

\begin{proposition} \label{pfw}
If in the pair $(A,B)$, $A\cap B=\emptyset$, with $A$ closed
decreasing and $B$ closed increasing, $A$ is functionally-preordered
closed and $B$ is functionally-preordered closed then $(A,B)$ is
functionally-preordered closed.
\end{proposition}

\begin{proof}
Let $A=g^{-1}(0)$ and $B=h^{-1}(1)$ with $g,h: E \to [0,1]$
continuous isotone functions. Let $\alpha: ([0,1]\times
[0,1])\backslash (0,1) \to [0,1]$ be a continuous function which is
isotone according to the product order on the square. Let $\alpha$
be also such that $\alpha^{-1}(0)=\{(0,y): y \in [0,1)\}$ and
$\alpha^{-1}(1)=\{((x,1): x\in (0,1]\}$. A possible example is
$\alpha=\frac{1+y}{2} x^{(1-y)/2}$; another example is
$\alpha=1/[1+(1-y)/x]$. The function $f=\alpha(g,h)$ is isotone,
continuous, and satisfies $f^{-1}(0)=g^{-1}(0)=A$,
$f^{-1}(1)=h^{-1}(0)=B$.
\end{proof}

The next result extends a known result for topological spaces
\cite[Theor. 16.8]{willard70}.

\begin{proposition} \label{hsr}
Every regularly preordered Lindel\"of space is a normally preordered
space.
\end{proposition}

\begin{proof}
Let $A$ and $B$ be closed disjoint sets which are respectively
decreasing and increasing. Since $A\cap B=\emptyset$, by preorder
regularity for each $x\in A$ there is an open decreasing set $U_x\ni
x$ such that $D(U_x)\cap B=\emptyset$. The collection $\{U_x\}$
covers $A$ and since the Lindel\"of property is hereditary with
respect to closed subspaces, there is a countable subcollection
$\{U_i\}$ with the same property. In the same way we find a
countable collection of open increasing sets $\{V_i\}$  which covers
$B$ and such that $I(V_i) \cap A=\emptyset$. Let us define the
sequence of open decreasing sets $W_1=U_1$,
$W_{n+1}=U_{n+1}\backslash [\cup_{i=1}^{n} I(V_i)]$, and the
sequence of open increasing sets $E_1=V_1\backslash D(U_1)$,
$E_{n}=V_{n}\backslash [\cup_{i=1}^{n}D(U_i)]$. The open disjoint
sets $U'=\cup_{n=1}^{\infty} W_n$ and $V'=\cup_{n=1}^{\infty} E_n$
are respectively decreasing and increasing and contain respectively
$A$ and $B$.
\end{proof}

\begin{theorem} \label{bha}
Every second countable regularly preordered space
$(E,\mathscr{T},\le)$ is a perfectly normally preordered space.
\end{theorem}

\begin{proof}

As second countability implies the Lindel\"of property, by
proposition \ref{hsr} $E$ is normally preordered. Let $A$ be a
closed decreasing set, we have only to prove that it is
functionally-preordered closed, the proof in the closed increasing
case being analogous. Suppose $A$ is open then $E\backslash A$ is
open, closed and increasing. Setting for $x \in A$, $f(x)=0$ and for
$x\in E\backslash A$, $f(x)=1$ we have finished. Therefore, we can
assume that $A$ is not open and hence that $A\ne E$.

Let $\mathcal{U}$ be a countable base of the topology of $E$. Let
$\mathcal{C} \subset \mathcal{U}$ be the subset whose elements,
denoted  $U_i$, $i \ge 1$, are such that $ A\cap I(U_i)=\emptyset$.
For every $U_k \in \mathcal{C}$, denote with $f_k: E \to [0,1]$ an
isotone continuous function which separates $A$ and $I(U_k)$, that
is, such that $f_{k}^{-1}(1)\supset I(U_{k})$ and
$f_{k}^{-1}(0)\supset A$ (see \cite[Theor. 1]{nachbin65}).

Let $x \in E\backslash A$. Applying  the preordered normality of the
space we can find  $\hat{U}$, open decreasing sets such that
$A\subset \hat{U} \subset D(\hat{U})  \subset E\backslash  i(x)$. As
$E\backslash D(\hat{U})$ is open there is some $O\in \mathcal{U}$
such that $x\in O \subset E\backslash D(\hat{U})$. Furthermore,
since $E\backslash \hat{U}\supset O$ is closed increasing,
$I(O)\subset E\backslash \hat{U}$ which implies $A\cap
I(O)=\emptyset$ and hence that there is some $U_k\in \mathcal{C}$
such that $O=U_k$. Thus there is a continuous isotone  function
$f_k$ such that $f_{k}^{-1}(1)\supset I(U_{k})$ and
$f_{k}^{-1}(0)\supset A$. In other words we have proved that for
each $x\in E\backslash A$ there is some $k$ such that $f_k(x)=1$.
Let us consider the function
\[
f=\sum_{k=1}^{\infty} \frac{1}{2^k} f_k .
\]
This function is clearly isotone, takes values in $[0,1]$ and it is
continuous because the series converges uniformly. Since for every
$k$, $f_{k}^{-1}(0)\supset A$, the same is true for $f$. Moreover,
note that if $x\in E\backslash A$ then there is some $k \ge 1$ such
that $f_k(x)>0$, which implies that $f(x)>0$ hence $f^{-1}(0)=A$.
\end{proof}

\begin{lemma} \label{nwx}
If a topological preordered space admits a countable continuous
isotone function representation, namely if there are continuous
isotone functions $g_k: E\to [0,1]$, $k\ge 1$ such that $ x\le y
\Leftrightarrow \forall k\ge 1, g_k(x)\le g_k(y)$, then there is a
countable continuous utility representation, namely  there are
continuous utility functions $f_k: E\to [0,1]$, $k\ge 1$ such that $
x\le y \Leftrightarrow \forall k\ge 1, f_k(x)\le f_k(y)$.
\end{lemma}

\begin{proof}
Let us consider the function
\[
g=\sum_{k=1}^{\infty} \frac{1}{2^k} g_k .
\]
This function is clearly isotone, takes values in $[0,1]$ and it is
continuous because the series converges uniformly. Furthermore, if
$x\le y$ and $y\nleq x$, then $g(x)<g(y)$ because all the $g_k$ are
isotone and there is some $g_{\bar{k}}$ for which
$g_{\bar{k}}(x)<g_{\bar{k}}(y)$. Indeed, if it were not so then for
every $k$, $g_{{k}}(y)\le g_{{k}}(x)$ which implies $y\le x$, a
contradiction. We conclude that $g$ is a continuous utility
function. Now, define for all $k,n\ge 1$,
$\tilde{g}_{kn}=(1-1/n)g_k+g/n$, so that $\tilde{g}_{kn}$ is a
continuous utility function. We have only to prove that if $x\nleq
y$ there are some $k,n,$ such that
$\tilde{g}_{kn}(x)>\tilde{g}_{kn}(y)$ but we know that there is a
$\tilde{k}$ such that $g_{\tilde{k}}(x)>g_{\tilde{k}}(y)$ (otherwise
for every $k$, $g_k(x)\le g_k(y)$ which implies $x\le y$, a
contradiction). Thus
$\tilde{g}_{\tilde{k}n}(x)=(1-1/n)g_{\tilde{k}}(x)+g(x)/n=
\tilde{g}_{\tilde{k}n}(y)+(1-1/n)(g_{\tilde{k}}(x)-g_{\tilde{k}}(y))+(g(x)-g(y))/n$
and thus the desired inequality holds for sufficiently large $n$.
\end{proof}

\begin{theorem} \label{bhb}
Every second countable regularly preordered space $E$ admits a
countable continuous  utility representation, that is, there is a
countable set $\{f_k, k\ge 1\}$ of continuous utility functions
$f_k: E\to [0,1]$  such that
\[
x\le y \Leftrightarrow \forall k\ge 1, f_k(x)\le f_k(y).
\]
\end{theorem}

\begin{proof}
Let $\mathcal{G}$ be the set of continuous isotone functions and for
every $g\in \mathcal{G}$ let $G_g=\{(x,y)\in E\times E: g(x) \le
g(y)\}$. The set $G_g$ is closed because of the continuity of $g$.
We already know that the preordered space $(E,\mathscr{T},\le)$ is a
normally preordered space (Theor. \ref{bha}), thus $\le$ is
represented by $\mathcal{G}$ (by \cite[Theor. 1]{nachbin65}), namely
$G(\le)=\bigcap_{g\in \mathcal{G}} G_g$. As $E$ is second countable,
$E\times E$ is second countable and hence hereditary Lindel\"of. As
a consequence, the intersection of an arbitrary family of closed
sets can be written as the intersection of a countable subfamily
$\mathcal{G}'\subset \mathcal{G}$. The desired conclusion follows
now from lemma \ref{nwx}.
\end{proof}

%
%

\section{Conclusions}

Every Hausdorff locally compact space is a completely regular space
\cite[Theor. 19.3]{willard70} and hence a regular space, and  every
second countable regular space ($T_3$-space) is metrizable
(Urysohn's theorem). These results allow us to improve the
separability properties of the space. Unfortunately, they have no
straightforward analogs in the general preordered case where a
quasi-uniformizable space need not be a regularly preordered space
\cite[Example 1]{kunzi94b}. One still would like to have some result
which improves the preorder separability properties of the space,
given suitable compactness or countability conditions on the
topology.
We have then followed  a different path proving that closed
preordered
 $k_\omega$-spaces are normally preordered. We have also proved
that second countable regularly preordered spaces are perfectly
normally preordered and that they admit a countable continuous
utility representation of the preorder.

\section*{Acknowledgments}
I thank A. Fathi for pointing our reference \cite{akin93}, and A.
Fathi, H.-P.A. K{\"u}nzi and  P. Pageault for stimulating
conversations on the subject of topological preordered spaces. I
also thank two anonymous referees for many useful suggestions. This
work has been partially supported by GNFM of INDAM and by FQXi.


\end{document}